\documentclass[reqno,12pt]{amsart}

\newtheorem{Theorem}{Theorem}[section]

\newtheorem{Corollary}[Theorem]{Corollary}
\newtheorem{Lemma}[Theorem]{Lemma}
\theoremstyle{remark}
\newtheorem{Remark}[Theorem]{Remark}
\numberwithin{equation}{section}

\allowdisplaybreaks

\begin{document}

\title[Curious $q$-expansions and beta integrals]
{Some curious $\boldsymbol q$-series expansions\\
and beta integral evaluations}

\author{George Gasper}
\address{Department of Mathematics, Northwestern University,
2033 Sheridan Road, Evanston, IL 60208-2730, USA}
\email{george@math.northwestern.edu}
\urladdr{http://www.math.northwestern.edu/{\textasciitilde}george}

\author[Michael Schlosser]{Michael Schlosser$^*$}
\address{Institut f\"ur Mathematik der Universit\"at Wien,
Nordbergstra{\ss}e 15, A-1090 Wien, Austria}
\email{schlosse@ap.univie.ac.at}
\urladdr{http://www.mat.univie.ac.at/{\textasciitilde}schlosse}

\thanks{$^*$The second author was fully supported by an APART
fellowship of the Austrian Academy of Sciences}
\date{March 29, 2004}

\dedicatory{Dedicated to Dick Askey
on the occasion of his 70th birthday}

\begin{abstract}
We deduce several curious $q$-series expansions by applying inverse
relations to certain identities for basic hypergeometric series.
After rewriting some of these expansions in terms of $q$-integrals,
we obtain, in the limit $q\to 1$, some curious beta-type integral
evaluations which appear to be new.
\end{abstract}

\maketitle

\section{Introduction}
Euler's beta integral evaluation (cf.~\cite[Eq.~(1.1.13)]{AAR})
\begin{equation} \label{betaf}
\int_0^1t^{\alpha-1}(1-t)^{\beta-1}dt=
\frac{\Gamma(\alpha)\Gamma(\beta)}{\Gamma(\alpha+\beta)},\qquad
\Re(a),\Re(b)>0,
\end{equation}
is one of the most important and prominent identities in
special functions. In Andrews, Askey and Roy's modern
treatise \cite{AAR}, the beta integral (and its various extensions)
runs like a thread through their whole exposition.

Concerning the evaluation of integrals of, say, elementary
functions, there is no general procedure which will find the
closed form evaluation if it exists. Being encountered with
some explicit integral (in this paper all integrals are definite),
if standard methods seem out of reach, it is usually wise
to consult one of the compound volumes listing
tables of integrals~\cite{E,PBM,RG}, hoping that the sought
evaluation could be found in there. However, this does not always
lead to success. In particular, several of the integral evaluations
obtained in this paper (specifically, Theorems~\ref{thbetat}
and \ref{thbetat2} and their specializations
\eqref{spec2}, \eqref{spec3}, \eqref{spec5})
are apparently not included as entries in the standard references
\cite{E,PBM,RG}. For instance, two special cases
($\alpha=\beta+1$, and $\alpha=\beta$, respectively) of one of our
main results (Theorem~\ref{thbetat}, which generalizes
\eqref{betaf}) are the following beta-type integral evaluations:
\begin{multline}\label{spec2}
\frac{\Gamma(\beta)\Gamma(\beta)}{2\,\Gamma(2\beta)}=
(c-(a+1)^2)\int_0^1
\frac{(c-a(a+t))^\beta\,(c-(a+1)(a+t))^{\beta-1}}
{(c-(a+t)^2)^{2\beta}}\\\times
t^{\beta}\,(1-t)^{\beta-1}\,dt
\end{multline}
and
\begin{multline}\label{spec3}
\frac{\Gamma(\beta)\Gamma(\beta)}{\Gamma(2\beta)}=
(c-(a+1)^2)\int_0^1
\frac{(c-a(a+t))^{\beta-1}\,(c-(a+1)(a+t))^{\beta-1}}
{(c-(a+t)^2)^{2\beta}}\\\times
(c-(a-t)(a+t))\, t^{\beta-1}\,(1-t)^{\beta-1}\,dt,
\end{multline}
where $\Re(\beta)>0$. (These evaluations and others
have been numerically verified using Mathematica.)

These integrals seem difficult to prove with standard methods,
such as expanding all factors in terms of powers of $t$
(by the binomial theorem) and integrating term-wise.
Applying this procedure to \eqref{spec2} yields a
five-fold sum that can be easily reduced to a four-fold sum,
but then one is apparently stuck.

In the sequel, we will develop some machinery
for proving our integral evaluations. First we derive, by
inverse relations, new $q$-series expansions. We then
rewrite these in terms of $q$-integrals. Finally, by letting $q\to 1$
we obtain the desired beta-type integral evaluations.

\section{Preliminaries}\label{secpre}

\subsection{Hypergeometric and basic hypergeometric series}

For a complex number $a$, define the {\em shifted factorial}
\begin{equation*}
(a)_0:=1,\qquad (a)_k:=a(a+1)\dots(a+k-1),
\end{equation*}
where $k$ is a positive integer. Let $r$ be a positive integer.
The {\em hypergeometric ${}_rF_{r-1}$ series}
with numerator parameters $a_1,\dots,a_r$, denominator parameters
$b_1,\dots,b_{r-1}$, and argument $z$ is defined by
\begin{equation*}
{}_rF_{r-1}\!\left[\begin{matrix}a_1,\dots,a_r\\
b_1,\dots,b_{r-1}\end{matrix};z\right]:=
\sum_{k\ge 0}\frac{(a_1)_k\dots(a_r)_k}
{k!\,(b_1)_k\dots(b_{r-1})_k}\,z^k{}.
\end{equation*}
The ${}_rF_{r-1}$ series terminates if one of the numerator parameters is
of the form $-n$ for a nonnegative integer $n$. If the series
does not terminate, it converges when $|z|<1$, and also when $|z|=1$ and
$\Re[b_1+b_2+\cdots+b_{r-1}-(a_1+a_2+\cdots+a_r)]>0$. See \cite{B,Sl}
for a classic texts on (ordinary) hypergeometric series.

Let $q$ (the ``base'') be a complex number such that $0<|q|<1$.
Define the {\em $q$-shifted factorial} by
\begin{equation*}
(a;q)_\infty:=\prod_{j\ge 0}(1-aq^j)\qquad\text{and}\qquad
(a;q)_k:=\frac{(a;q)_\infty}{(aq^k;q)_\infty}
\end{equation*}
for integer $k$.
The {\em basic hypergeometric ${}_r\phi_{r-1}$ series}
with numerator parameters $a_1,\dots,a_r$, denominator parameters
$b_1,\dots,b_{r-1}$, base $q$, and argument $z$ is defined by
\begin{equation*}
{}_r\phi_{r-1}\!\left[\begin{matrix}a_1,\dots,a_r\\
b_1,\dots,b_{r-1}\end{matrix};q,z\right]:=
\sum_{k\ge 0}\frac{(a_1;q)_k\dots(a_r;q)_k}
{(q;q)_k(b_1;q)_k\dots(b_{r-1};q)_k}\,z^k.
\end{equation*}
The ${}_r\phi_{r-1}$ series terminates if one of the numerator parameters
is of the form $q^{-n}$ for a nonnegative integer $n$. If the series
does not terminate, it converges when $|z|<1$.
For a thorough exposition on basic hypergeometric series
(or, synonymously, {\em $q$-hypergeometric series}),
including a list of several selected summation and transformation formulas,
we refer the reader to \cite{GR}.

We list three specific identities which we will utilize in this paper.

\begin{Lemma}[$q$-Kummer summation (cf.\ {\cite[Eq.\ (II.9)]{GR}})]
\label{qkummer}
\begin{equation*}
{}_2\phi_1\!\left[\begin{matrix}a,b\\
aq/b\end{matrix};q,-\frac qb\right]=
\frac{(-q;q)_\infty\,(aq;q^2)_\infty(aq^2/b^2;q^2)_\infty}
{(-q/b;q)_\infty(aq/b;q)_\infty},
\end{equation*}
provided $|q/b|<1$.
\end{Lemma}
\begin{proof}
One may simply specialize Rogers' nonterminating
very-well-poised ${}_6\phi_5$ summation
(cf.\ \cite[Eq.~(II.20)]{GR})
\begin{multline*}
{}_6\phi_5\!\left[\begin{matrix}a,\,q\sqrt{a},-q\sqrt{a},b,c,d\\
\sqrt{a},-\sqrt{a},aq/b,aq/c,aq/d\end{matrix};q,
\frac{aq}{bcd}\right]\\
=\frac {(aq;q)_{\infty}(aq/bc;q)_{\infty}(aq/bd;q)_{\infty}
(aq/cd;q)_{\infty}}
{(aq/b;q)_{\infty}(aq/c;q)_{\infty}(aq/d;q)_{\infty}(aq/bcd;q)_{\infty}},
\end{multline*}
where $|aq/bcd|<1$, by setting $c=\sqrt{a}$ and $d=-\sqrt{a}$,
hereby ``cancelling off'' the very-well-poised term.
\end{proof}

For a simple derivation of the following transformation from
the $q$-binomial theorem \eqref{qbin}, see \cite[\S~1.4]{GR}.

\begin{Lemma}[Second iterate of Heine's transformation
(cf.\ {\cite[Eq.\ (III.2)]{GR}})]\label{heine}
\begin{equation*}
{}_2\phi_1\!\left[\begin{matrix}a,b\\
c\end{matrix};q,z\right]=
\frac{(c/b;q)_\infty\,(bz;q)_\infty}
{(c;q)_\infty(z;q)_\infty}\,
{}_2\phi_1\!\left[\begin{matrix}abz/c,b\\
bz\end{matrix};q,\frac cb\right],
\end{equation*}
provided $|z|,|c/b|<1$.
\end{Lemma}

\begin{Lemma}[An ($m+1$)-term ${}_3\phi_2$ summation]\label{lem2}
Let $m$ be a nonnegative integer. Then
\begin{equation}\label{lem2id}
{}_3\phi_2\!\left[\begin{matrix}a,b,dq^m\\
c,d\end{matrix};q,\frac{cq^{-m}}{ab}\right]=
\frac{(c/a;q)_\infty(c/b;q)_\infty}
{(c;q)_\infty(c/ab;q)_\infty}\,
{}_3\phi_2\!\left[\begin{matrix}a,b,q^{-m}\\
abq/c,d\end{matrix};q,q\right],
\end{equation}
provided $|cq^{-m}/ab|<1$.
\end{Lemma}
\begin{proof}
This can be obtained from \cite[Eq.~(III.34)]{GR}, i.e.\ the
three-term transformation
\begin{multline*}
{}_3\phi_2\!\left[\begin{matrix}a,b,c\\
d,e\end{matrix};q,\frac{de}{abc}\right]=
\frac{(e/b;q)_\infty(e/c;q)_\infty}
{(e;q)_\infty(e/bc;q)_\infty}\,
{}_3\phi_2\!\left[\begin{matrix}d/a,b,c\\
d,bcq/e\end{matrix};q,q\right]\\
+\frac{(d/a;q)_\infty(b;q)_\infty(c;q)_\infty(de/bc;q)_\infty}
{(d;q)_\infty(e;q)_\infty(bc/e;q)_\infty(de/abc;q)_\infty}\,
{}_3\phi_2\!\left[\begin{matrix}e/b,e/c,de/abc\\
de/bc,eq/bc\end{matrix};q,q\right],
\end{multline*}
where $|de/abc|<1$, by first letting $a\to dq^m$, by which the
coefficient of the second ${}_3\phi_2$ on the right-hand vanishes,
and then suitably relabeling the parameters.
\end{proof}
A more direct proof of Lemma~\ref{lem2} proceeds by induction on $m$,
using the $q$-Gau{\ss} summation \eqref{qgauss} in the inductive basis,
and the simple identity
\begin{equation*}
\frac{1-dq^{m+k}}{1-dq^m}=q^k+\frac{1-q^k}{1-dq^m}
\end{equation*}
in the inductive step. The details are left to the reader.

\begin{Remark}
We view \eqref{lem2id} as a {\em summation} (versus a
{\em transformation}) since the left-hand side is a nonterminating
sum and the right hand-side contains a finite number of terms.
For $m=0$ Lemma~\ref{lem2} reduces to the classical
$q$-Gau{\ss} summation (cf.\ {\cite[Eq.\ (II.8)]{GR}})
\begin{equation}\label{qgauss}
{}_2\phi_1\!\left[\begin{matrix}a,b\\
c\end{matrix};q,\frac c{ab}\right]=
\frac{(c/a;q)_\infty(c/b;q)_\infty}
{(c;q)_\infty(c/ab;q)_\infty},
\end{equation}
whereas for $m=1$ it reduces to
\begin{equation*}
{}_3\phi_2\!\left[\begin{matrix}a,b,dq\\
c,d\end{matrix};q,\frac c{abq}\right]=
\left(1-\frac{(1-a)(1-b)}{(1-abq/c)(1-d)}\right)
\frac{(c/a;q)_\infty(c/b;q)_\infty}
{(c;q)_\infty(c/ab;q)_\infty}.
\end{equation*}
\end{Remark}

\subsection{Inverse relations}

Let $\mathbb Z$ denote the set of integers
and  $F=(f_{nk})_{n,k\in\mathbb Z}$ be an infinite lower-triangular
matrix; i.e.\ $f_{nk}=0$ unless $n\ge k$.
The matrix $G=(g_{kl})_{k,l\in\mathbb Z}$ is said
to be the {\em inverse matrix} of $F$ if and only if
\begin{equation*}
\sum_{l\le k\le n} f_{nk}g_{kl}=\delta_{nl}
\end{equation*}
for all $n,l\in\mathbb Z$, where $\delta_{nl}$ is the
usual Kronecker delta.

The method of applying {\em inverse relations} \cite{Ri}
is a well-known technique for proving identities, or for
producing new ones from given ones.

If $(f_{nk})_{n,k\in\mathbb Z}$ and $(g_{kl})_{k,l\in\mathbb Z}$ are
lower-triangular matrices that are inverses of each other,
then
\begin{subequations}\label{rotinv}
\begin{equation}\label{rotinvf}
\sum_{n\ge k}f_{nk}a_n=b_k
\end{equation}
{\em if and only if}
\begin{equation}\label{rotinvg}
\sum_{k\ge l}
g_{kl}b_k=a_l,
\end{equation}
\end{subequations}
subject to suitable convergence conditions.
For some applications of \eqref{rotinv} see e.g.\ \cite{Kr,Ri,Sc}.

Note that in the literature it is actually more common to consider
the following inverse relations involving finite sums,
\begin{equation}\label{inv}
\sum_{k=0}^nf_{nk}a_k=b_n\qquad\text{\em if and only if}\qquad
\sum_{l=0}^kg_{kl}b_l=a_k.
\end{equation}

It is clear that in order to apply \eqref{rotinv} (or \eqref{inv})
effectively, one should have some explicit matrix inversion at hand.

\begin{Lemma}[Krattenthaler~\cite{Kr}]\label{kmi}
Let $(a_j)_{j\in\mathbb Z}$, $(c_j)_{j\in\mathbb Z}$
be arbitrary sequences and $d$ an arbitrary indeterminate. Then
the infinite matrices
$(f_{nk})_{n,k\in\mathbb Z}$ and
$(g_{kl})_{k,l\in\mathbb Z}$ are inverses of each other, where
\begin{subequations}
\begin{equation*}
f_{nk}=\frac{\prod_{j=k}^{n-1}(a_j-d/c_k)(a_j-c_k)}
{\prod_{j=k+1}^n(c_j-d/c_k)(c_j-c_k)},
\end{equation*}
\begin{equation*}
g_{kl}=\frac{(a_lc_l-d)(a_l-c_l)}
{(a_kc_k-d)(a_k-c_k)}
\frac{\prod_{j=l+1}^k(a_j-d/c_k)(a_j-c_k)}
{\prod_{j=l}^{k-1}(c_j-d/c_k)(c_j-c_k)}.
\end{equation*}
\end{subequations}
\end{Lemma}

Krattenthaler's matrix inverse is very general as
it contains a vast number of other known explicit infinite matrix
inversions. Several of its useful special cases are of (basic)
hypergeometric type. The following special case of Lemma~\ref{kmi}
is exceptional in the sense that although it involves powers of $q$,
it is {\em not} to be considered a $q$-hypergeometric inversion.
(More precisely, the following special case serves as a bridge
between $q$-hypergeometric and certain non-$q$-hypergeometric
identities. For some other such matrix inverses, see \cite{Sc}.)

\begin{Corollary}[MS {\cite[Eqs.~(7.18)/(7.19)]{Sc}}]\label{cor1}
Let
\begin{subequations}
\begin{equation*}
f_{n k}=\frac {(1/b;q)_{n-k}\,
\big(\frac {(a+bq^k)q^k}{c-a(a+bq^k)};q\big)_{n-k}}
{(q;q)_{n-k}\,\big(\frac {(a+bq^k)bq^{k+1}}{c-a(a+bq^k)};q\big)_{n-k}},
\end{equation*}
\begin{equation*}
g_{k l}=(-1)^{k-l}\,q^{\binom{k-l}2}
\frac{(c-(a+bq^l)(a+q^l))}{(c-(a+bq^k)(a+q^k))}
\frac{(q^{l-k+1}/b;q)_{k-l}\,
\big(\frac {(a+bq^k)q^{l+1}}{c-a(a+bq^k)};q\big)_{k-l}}
{(q;q)_{k-l}\,\big(\frac {(a+bq^k)bq^l}{c-a(a+bq^k)};q\big)_{k-l}}.
\end{equation*}
\end{subequations}
Then the infinite matrices $(f_{nk})_{n,k\in\mathbb Z}$ and
$(g_{kl})_{k,l\in\mathbb Z}$ are inverses of each other.
\end{Corollary}
\begin{proof}
In Lemma~\ref{kmi} set $a_j\mapsto a+q^j$, $c_j\mapsto a+bq^j$
($j\in\mathbb Z$), and $d\mapsto c$,
and perform some elementary manipulations.
\end{proof}

\section{Some curious $q$-series expansions}\label{secqexp}

Corollary~\ref{cor1} was utilized in \cite[Th.~7.16]{Sc} to deduce
from the classical $q$-Gau{\ss} summation a specific
$q$-series expansion, the latter itself not belonging to the
hierarchy of basic hypergeometric series.
In particular, the following identity was obtained:
\begin{multline}\label{ntnewq0gl}
\frac {(b^2q;q)_{\infty}}{(bq;q)_{\infty}}=
\sum_{k=0}^{\infty}\frac {(c-(a+1)(a+b))}{(c-(a+1)(a+bq^k))}\,
\frac {(c-(a+bq^k)^2)}{(c-(a+b)(a+bq^k))}\\\times
\frac {(b;q)_k\,\big(\frac {(a+bq^k)}{c-a(a+bq^k)};q\big)_k\,
\big(\frac {(a+bq^k)b^2q^{k+1}}{c-a(a+bq^k)};q\big)_{\infty}}
{(q;q)_k\,\big(\frac {(a+bq^k)bq}{c-a(a+bq^k)};q\big)_{\infty}}\,
(bq)^k,
\end{multline}
where $|bq|<1$. We find this to be quite a curious expansion.
For $c=0$ it reduces to a special case of the $q$-Gau{\ss}
summation \eqref{qgauss}. On the other hand, for $a=0$ it reduces to
\begin{equation*}
\frac {(b^2q;q)_{\infty}(b^2q/c;q)_{\infty}}
{(bq;q)_{\infty}(b^3q/c;q)_{\infty}}=
\sum_{k=0}^{\infty}\frac{(1-b^2q^{2k}/c)}{(1-b^2/c)}
\frac{(b^2/c;q)_k(b;q)_k\,(b/c;q)_{2k}}
{(q;q)_k(bq/c;q)_k\,(b^3q/c;q)_{2k}}\,(bq)^k,
\end{equation*}
a particular very-well-poised ${}_8\phi_7$ summation,
which is equivalent to the $n\to\infty$ special case of
the terminating very-well-poised ${}_{10}\phi_9$ summation,
\begin{multline}\label{212}
{}_{10}\phi_9\!\left[\begin{matrix}a,\,q\sqrt{a},-q\sqrt{a},
\sqrt{b},-\sqrt{b},\sqrt{bq},-\sqrt{bq},a/b,a^2q^{n+1}/b,q^{-n}\\
\sqrt{a},-\sqrt{a},aq/\sqrt{b},-aq/\sqrt{b},a\sqrt{q/b},-a\sqrt{q/b},
bq,bq^{-n}/a,aq^{n+1}\end{matrix};q,q\right]\\
=\frac{(aq;q)_n\,(a^2q/b^2;q)_n}{(aq/b;q)_n\,(a^2q/b;q)_n},
\end{multline}
given in \cite[Ex.~2.12]{GR}. This ${}_{10}\phi_9$ summation
itself follows immediately from taking the limit $d\to 1$ in
\cite[Eq.~(2.8.3)]{GR} which is Bailey's~\cite[p.~431]{Ba1},
\cite{Ba2} transformation of a terminating, balanced,
nearly-poised of the second kind ${}_5\phi_4$ series into
a multiple of a particular terminating, balanced,
very-well-poised ${}_{12}\phi_{11}$ series.
(See \cite{GR} for the terminology.)
This latter ${}_5\phi_4\leftrightarrow{}_{12}\phi_{11}$
transformation is a consequence of the
``WP-Bailey lemma'', cf.\ \cite[Eq.~(2.8.2)]{GR} and
\cite[\S\S~6 and 7]{An}.

The matrix inversion in Corollary~\ref{cor1} was also applied
to both the classical $q$-Pfaff--Saalsch\"utz summation and the
$2$-balanced ${}_3\phi_2$ summation, to derive
two non-$q$-hypergeometric terminating summations,
see \cite[Ths.~7.34 and 7.38]{Sc}. For illustration
(and to correct some misprints which appeared in the printed
version of \cite{Sc}), we reproduce the first one of these,
specifically, (7.35) of \cite[Th.~7.34]{Sc}:
\begin{multline}\label{newq01gl}
\frac {(c^2q;q)_n} {(cq;q)_n}=
\sum_{k=0}^n\frac{(b+(a-c)(a-1))}{(b+(a-c)(a-q^{-k}))}
\frac{(b+(a-q^{-k})^2)}{(b+(a-1)(a-q^{-k}))}\\\times
\frac{(q^{-n};q)_k\,(c;q)_k\,
\big(\frac{b+a(a-q^{-k})}{c(a-q^{-k})};q\big)_k}
{(q;q)_k\,(q^{-n}/c;q)_k\,
\big(cq\frac{b+a(a-q^{-k})}{(a-q^{-k})};q\big)_k}\,
\frac{\big(cq\frac{b+a(a-q^{-k})}{(a-q^{-k})};q\big)_n}
{\big(q\frac{b+a(a-q^{-k})}{(a-q^{-k})};q\big)_n}\,q^k.
\end{multline}
The $n\to\infty$ case of \eqref{newq01gl} is equivalent to
\eqref{ntnewq0gl}. For $a=0$, \eqref{newq01gl} is equivalent
to \eqref{212}.

For other terminating and nonterminating summations that were
derived via inverse relations from
classical ordinary and basic hypergeometric summations and
do {\em not} belong to the hierarchy of (basic)
hypergeometric series, such as identities of ($q$-)Abel,
($q$-)Rothe, or of the above type (as in \eqref{ntnewq0gl} and
\eqref{newq01gl}), see \cite{Sc}.

We commence with a new application of Corollary~\ref{cor1}, which
was missed in \cite[\S~7]{Sc}.

\begin{Theorem}\label{th1}
Let $a$, $b$, and $c$ be indeterminate.
Then
\begin{multline}\label{th1id}
\frac {(-bq;q)_{\infty}}{(-q;q)_{\infty}}=
\sum_{k=0}^{\infty}\frac {(c-(a+1)(a+b))}{(c-(a+1)(a+bq^k))}\,
\frac {(c-(a+bq^k)^2)}{(c-(a+b)(a+bq^k))}\\\times
\frac {(b;q)_k\,\big(\frac {(a+bq^k)}{c-a(a+bq^k)};q\big)_\infty}
{(q;q)_k\,\big(\frac {(a+bq^k)bq}{c-a(a+bq^k)};q\big)_{\infty}}\,
\frac{\big(\frac {(a+bq^k)b^2q^{k+2}}{c-a(a+bq^k)};q^2\big)_{\infty}}
{\big(\frac {(a+bq^k)q^k}{c-a(a+bq^k)};q^2\big)_\infty}\,
(-1)^k\,q^k.
\end{multline}
\end{Theorem}

For $a=0$, \eqref{th1id} reduces to a special case of Rogers'
very-well-poised ${}_6\phi_5$ summation (cf.\ \cite[Eq.~(II.20)]{GR}).
On the other hand, if $c=0$, we obtain (with $a\mapsto -1/a$)
\begin{equation}\label{quid}
\frac{(-bq;q)_\infty\,(abq;q)_\infty}
{(-q;q)_\infty\,(a;q)_\infty}=\sum_{k=0}^{\infty}
\frac{(b;q)_k}{(q;q)_k}\frac{(ab^2q^{2+k};q^2)_\infty}
   {(aq^k;q^2)_\infty}\,(-1)^k\,q^k,
\end{equation}
which we could not find in this explicit form in the literature.
Nevertheless, it is not difficult to find a conventional proof.
Splitting the sum on the right hand side in two parts depending
on the parity of $k$, \eqref{quid} becomes
\begin{multline*}\frac{(-bq;q)_\infty\,(abq;q)_\infty}
{(-q;q)_\infty\,(a;q)_\infty}=
\frac{(ab^2q^2;q^2)_\infty}{(a;q^2)_\infty}\,
{}_3\phi_2\!\left[\begin{matrix}b,bq,a\\
q,ab^2q^2\end{matrix};q^2,q^2\right]\\-
q\frac{(1-b)}{(1-q)}\frac{(ab^2q^3;q^2)_\infty}{(aq;q^2)_\infty}\,
{}_3\phi_2\!\left[\begin{matrix}bq,bq^2,aq\\
q^3,ab^2q^3\end{matrix};q^2,q^2\right].
\end{multline*}
Now, this is just a special case of the nonterminating balanced
${}_3\phi_2$ summation \cite[Eq.~(II.24)]{GR}.

Other {\em quadratic} identities similar to \eqref{quid}
(with infinite products in the summand) have been derived in
\cite[Th.~4.2, Cors.~5.4, 5.5 and 5.6]{Sc2}.

\begin{proof}[Proof of Theorem~\ref{th1}]
Let the inverse matrices $(f_{nk})_{n,k\in\mathbb Z}$ and
$(g_{kl})_{k,l\in\mathbb Z}$ be defined as in Corollary~\ref{cor1}.
Then \eqref{rotinvf} holds for
\begin{equation*}
a_n=(-bq)^n\quad\text{and}\quad
b_k=(-bq)^k\frac {(-q;q)_{\infty}
\big(\frac {(a+bq^k)q^{k+1}}{c-a(a+bq^k)};q^2\big)_{\infty}\!
\big(\frac {(a+bq^k)b^2q^{k+2}}{c-a(a+bq^k)};q^2\big)_{\infty}}
{(-bq;q)_{\infty}
\big(\frac {(a+bq^k)bq^{k+1}}{c-a(a+bq^k)};q\big)_{\infty}}
\end{equation*}
by Lemma~\ref{qkummer}.
This implies the inverse relation \eqref{rotinvg},
with the above values of $a_n$ and $b_k$.
After performing the shift $k\mapsto k+l$, and
the substitutions $a\mapsto aq^l$, $c\mapsto cq^{2l}$, we get
rid of $l$ and eventually obtain \eqref{th1id}.
\end{proof}

Next, we present two generalizations of \eqref{ntnewq0gl}.

\begin{Theorem}\label{th}
Let $a$, $b$, and $c$ be indeterminate. Then
\begin{multline}\label{thid}
\frac {(z;q)_\infty}{(z/b;q)_\infty}=
\sum_{k=0}^{\infty}\frac {(c-(a+1)(a+b))}{(c-(a+1)(a+bq^k))}\,
\frac {(c-(a+bq^k)^2)}{(c-(a+b)(a+bq^k))}\\\times
{}_2\phi_1\!\left[\begin{matrix}1/b,z/b^2q\\
z/b\end{matrix};q,\frac {(a+bq^k)b^2q^{k+1}}{c-a(a+bq^k)}\right]\\\times
\frac {(b;q)_k\,\big(\frac {(a+bq^k)}{c-a(a+bq^k)};q\big)_k\,
\big(\frac {(a+bq^k)b^2q^{k+1}}{c-a(a+bq^k)};q\big)_{\infty}}
{(q;q)_k\,\big(\frac {(a+bq^k)bq}{c-a(a+bq^k)};q\big)_{\infty}}\,
(z/b)^k,
\end{multline}
provided $|z/b|<1$.
\end{Theorem}

Clearly, \eqref{thid} reduces to \eqref{ntnewq0gl} when $z=b^2q$.
At first glance, it seems that \eqref{thid} is not at all related
to \eqref{th1id} which also contains the base $q^2$. Notwithstanding,
\eqref{thid} is indeed more general than \eqref{th1id} and
reduces to the latter for $z=-bq$. In this case the ${}_2\phi_1$
appearing in the summand of \eqref{thid} becomes a ${}_1\phi_0$
with base $q^2$ (using $(1/b;q)_k(-1/b;q)_k=(1/b^2;q^2)_k$, etc.)
which can be summed by virtue of \eqref{qbin}.

\begin{proof}[Proof of Theorem~\ref{th}]
Let the inverse matrices $(f_{nk})_{n,k\in\mathbb Z}$ and
$(g_{kl})_{k,l\in\mathbb Z}$ be defined as in Corollary~\ref{cor1}.
Then \eqref{rotinvf} holds for $a_n=z^n$ and
\begin{equation*}
b_k=z^k\,\frac{(z/b;q)_\infty\,
\big(\frac{(a+bq^k)b^2q^{k+1}}{c-a(a+bq^k)};q\big)_{\infty}}
{(z;q)_\infty\,\big(\frac{(a+bq^k)bq^{k+1}}{c-a(a+bq^k)};q\big)_{\infty}}\,
{}_2\phi_1\!\left[\begin{matrix}1/b,z/b^2q\\
z/b\end{matrix};q,\frac{(a+bq^k)b^2q^{k+1}}{c-a(a+bq^k)}\right]
\end{equation*}
by Lemma~\ref{heine}.
This implies the inverse relation \eqref{rotinvg},
with the above values of $a_n$ and $b_k$.
After performing the shift $k\mapsto k+l$, and
the substitutions $a\mapsto aq^l$, $c\mapsto cq^{2l}$,
we get rid of $l$ and eventually obtain \eqref{thid}.
\end{proof}

\begin{Theorem}\label{th2}
Let $a$, $b$, and $c$ be indeterminate, and
let $m$ be a nonnegative integer. Then
\begin{multline}\label{th2id}
\frac {(b^2q;q)_{\infty}}{(bq;q)_{\infty}}=
\sum_{k=0}^{\infty}\frac {(c-(a+1)(a+b))}{(c-(a+1)(a+bq^k))}\,
\frac {(c-(a+bq^k)^2)}{(c-(a+b)(a+bq^k))}\\\times
{}_3\phi_2\!\left[\begin{matrix}1/b,
\frac{(a+bq^k)q^k}{c-a(a+bq^k)},q^{-m}\\[.4em]
1/b^2,eq^k\end{matrix};q,q\right]\,
\frac{(eq^m;q)_k}{(e;q)_k}\\\times
\frac {(b;q)_k\,\big(\frac {(a+bq^k)}{c-a(a+bq^k)};q\big)_k\,
\big(\frac {(a+bq^k)b^2q^{k+1}}{c-a(a+bq^k)};q\big)_{\infty}}
{(q;q)_k\,\big(\frac {(a+bq^k)bq}{c-a(a+bq^k)};q\big)_{\infty}}\,
(bq^{1-m})^k,
\end{multline}
provided $|bq^{1-m}|<1$.
\end{Theorem}

Clearly, \eqref{th2id} reduces to \eqref{ntnewq0gl} when $m=0$,
or when $e\to\infty$.

\begin{proof}[Proof of Theorem~\ref{th2}]
Let the inverse matrices $(f_{nk})_{n,k\in\mathbb Z}$ and
$(g_{kl})_{k,l\in\mathbb Z}$ be defined as in Corollary~\ref{cor1}.
Then \eqref{rotinvf} holds for
\begin{multline*}
a_n=\frac{(eq^m;q)_n}{(e;q)_n}\,(b^2q^{1-m})^n\qquad\text{and}\qquad
b_k={}_3\phi_2\!\left[\begin{matrix}1/b,
\frac{(a+bq^k)q^k}{c-a(a+bq^k)},q^{-m}\\[.4em]
1/b^2,eq^k\end{matrix};q,q\right]\\\times
\frac{(eq^m;q)_k}{(e;q)_k}\,
\frac {(bq;q)_{\infty}\,
\big(\frac {(a+bq^k)b^2q^{k+1}}{c-a(a+bq^k)};q\big)_{\infty}}
{(b^2q;q)_{\infty}\,
\big(\frac {(a+bq^k)bq^{k+1}}{c-a(a+bq^k)};q\big)_{\infty}}\,
(b^2q^{1-m})^k
\end{multline*}
by Lemma~\ref{lem2}.
This implies the inverse relation \eqref{rotinvg},
with the above values of $a_n$ and $b_k$.
After performing the shift $k\mapsto k+l$, and
the substitutions $a\mapsto aq^l$, $c\mapsto cq^{2l}$,
$e\mapsto eq^{-l}$, we get
rid of $l$ and eventually obtain \eqref{th2id}.
\end{proof}

\section{$q$-Integrals}\label{secqint}

In the following we restrict ourselves to {\em real} $q$ with $0<q<1$.

Thomae~\cite{Th} introduced the $q$-integral defined by
\begin{equation}\label{qint}
\int_0^1 f(t)d_qt=(1-q)\sum_{k=0}^{\infty}f(q^k)q^k.
\end{equation}
Later Jackson~\cite{Ja} gave a more general $q$-integral which
however we do not need here.

By considering the Riemann sum for a continuous function $f$
over the closed interval $[0,1]$, partitioned by the points $q^k$,
$k\ge 0$, one easily sees that
\begin{equation*}
\lim_{q\to 1^-}\int_0^1 f(t)d_qt=\int_0^1 f(t)dt.
\end{equation*}

It is well known that many identities for $q$-series can be written
in terms of $q$-integrals, which then may be specialized
(as $q\to 1$) to ordinary integrals. For instance,
the $q$-binomial theorem (cf.\ \cite[Eq.~(II.3)]{GR})
\begin{equation}\label{qbin}
\sum_{k=0}^{\infty}\frac{(a;q)_k}{(q;q)_k}z^k=
\frac{(az;q)_{\infty}}{(z;q)_{\infty}},\qquad |z|<1,
\end{equation}
can be written, when $a\mapsto q^{\beta}$ and $z\mapsto q^{\alpha}$, as
\begin{equation}\label{qbeta}
\int_0^1\frac{(qt;q)_\infty}{(q^{\beta}t;q)_\infty}t^{\alpha-1}d_qt=
\frac{\Gamma_q(\alpha)\Gamma_q(\beta)}{\Gamma_q(\alpha+\beta)},
\end{equation}
where
\begin{equation}\label{gammaq}
\Gamma_q(x):=(1-q)^{1-x}\frac{(q;q)_\infty}{(q^x;q)_\infty}
\end{equation}
is the $q$-gamma function, introduced by Thomae~\cite{Th},
see also \cite[\S~10.3]{AAR} and \cite[\S~1.11]{GR}. In fact,
\eqref{qbeta} is a $q$-extension of the beta integral
evaluation \eqref{betaf}.

Since the expansion of Proposition~\ref{th1} involves an alternating
series, it makes no sense to rewrite it as a $q$-integral;
the limit $q\to 1$ would never produce a convergent integral.
However, we can reasonably rewrite the expansions in Theorems~\ref{th}
and \ref{th2} in terms of $q$-integrals. These will then be utilized
in Section~\ref{secint} to obtain new beta-type integral evaluations.

Starting with \eqref{thid}, if we replace $z$ by $q^{\alpha+\beta}$,
$b$ by $q^{\beta}$, and multiply both sides of the identity by
\begin{equation}\label{factor}
(1-q)\frac{(q;q)_\infty}{(q^{\beta};q)_\infty},
\end{equation}
we obtain the following generalization of the $q$-beta
integral evaluation:
\begin{multline}\label{qbetat}
\frac{\Gamma_q(\alpha)\Gamma_q(\beta)}{\Gamma_q(\alpha+\beta)}=
\int_0^1 \frac {(c-(a+1)(a+q^\beta))}{(c-(a+1)(a+q^\beta t))}
\frac {(c-(a+q^\beta t)^2)}{(c-(a+q^\beta)(a+q^\beta t))}\\\times
{}_2\phi_1\!\left[\begin{matrix}q^{\alpha-\beta-1},q^{-\beta}\\
q^{\alpha}\end{matrix};q,
\frac{(a+q^\beta t)q^{2\beta+1}t}{c-a(a+q^\beta t)}\right]\\\times
\frac {(qt;q)_\infty\,
\big(\frac{(a+q^\beta t)}{c-a(a+q^\beta t)};q\big)_\infty\,
\big(\frac{(a+q^\beta t)q^{2\beta+1}t}{c-a(a+q^\beta t)};q\big)_{\infty}}
{(q^\beta t;q)_\infty\,
\big(\frac {(a+q^\beta t)t}{c-a(a+q^\beta t)};q\big)_\infty\,
\big(\frac {(a+q^\beta t)q^{\beta+1}}{c-a(a+q^\beta t)};q\big)_{\infty}}\,
t^{\alpha-1} d_qt.
\end{multline}
Clearly, \eqref{qbetat} reduces to \eqref{qbeta} when either
$c\to\infty$ or $a\to\infty$.

Similarly, starting with \eqref{th2id},
if we replace $b$ by $q^{\beta}$ and multiply
both sides of the identity by \eqref{factor},
we obtain the following $q$-beta-type integral evaluation:
\begin{multline}\label{qbetat2}
\frac{\Gamma_q(\beta+1)\Gamma_q(\beta)}{\Gamma_q(2\beta+1)}=
\int_0^1 \frac {(c-(a+1)(a+q^\beta))}{(c-(a+1)(a+q^\beta t))}
\frac {(c-(a+q^\beta t)^2)}{(c-(a+q^\beta)(a+q^\beta t))}\\\times
{}_3\phi_2\!\left[\begin{matrix}q^{-\beta},
\frac{(a+q^\beta t)t}{c-a(a+q^\beta t)},q^{-m}\\[.4em]
q^{-2\beta},et\end{matrix};q,q\right]\,
\frac{(eq^m;q)_\infty}{(e;q)_\infty}\,
\frac{(et;q)_\infty}{(eq^mt;q)_\infty}\\\times
\frac {(qt;q)_\infty\,
\big(\frac {(a+q^\beta t)}{c-a(a+q^\beta t)};q\big)_\infty\,
\big(\frac {(a+q^\beta t)q^{2\beta+1}t}{c-a(a+q^\beta t)};q\big)_{\infty}}
{(q^\beta t;q)_\infty\,
\big(\frac {(a+q^\beta t)t}{c-a(a+q^\beta t)};q\big)_\infty\,
\big(\frac {(a+q^\beta t)q^{\beta+1}}{c-a(a+q^\beta t)};q\big)_{\infty}}\,
t^{\beta-m} d_qt.
\end{multline}
This formula does not really extend \eqref{qbeta} as there is
only one ``exponent parameter'', $\beta$. However, for small $m$ the
${}_3\phi_2$ appearing in the integrand can be expanded in explicit
terms; the integrand thus has nearly ``closed form''.
In particular, the ($e=0$ and) $m=0$ case of \eqref{qbetat2} is equal to
the $\alpha=\beta+1$ special case of \eqref{qbetat}. More generally,
for $e=0$ equation \eqref{qbetat2} is equivalent to the $\alpha=\beta+1-m$
special case of \eqref{qbetat}, due to the
transformation (cf.\ \cite[Eq.~(III.7)]{GR})
\begin{equation*}
{}_2\phi_1\!\left[\begin{matrix}q^{-m},b\\c\end{matrix};q,z\right]=
\frac{(c/b;q)_m}{(c;q)_m}\,
{}_3\phi_2\!\left[\begin{matrix}q^{-m},b,bzq^{-m}/c\\
bq^{1-m}/c,0\end{matrix};q,q\right].
\end{equation*}

\section{Curious beta-type integrals}\label{secint}

Observe that $\lim_{q\to 1^-}\Gamma_q(x)=\Gamma(x)$
(see \cite[(1.10.3)]{GR}) and
\begin{equation*}
\lim_{q\to 1^-}\frac{(q^{\alpha}u;q)_\infty}{(u;q)_\infty}=
(1-u)^{-\alpha}
\end{equation*}
for constant $u$ (with $|u|<1$), due to \eqref{qbin} and its
$q\to 1$ limit, the ordinary binomial theorem.

We thus immediately deduce, as consequences of our
$q$-integral evaluations from Section~\ref{secqint},
some beta-type integral evaluations.
Throughout it is implicitly assumed that the integrals are well defined,
in particular that the parameters are chosen such that no poles
occur on the path of integration $t\in[0,1]$ and the integrals
converge.

The first beta integral evaluation is obtained from letting
$q\to 1$ in \eqref{qbetat}.
\begin{Theorem}\label{thbetat}
Let $\Re(\beta),\Re(\alpha)>0$. Then
\begin{multline}\label{betat}
\frac{\Gamma(\alpha)\Gamma(\beta)}{\Gamma(\alpha+\beta)}=
(c-(a+1)^2)\int_0^1
\frac{(c-a(a+t))^\beta\,(c-(a+1)(a+t))^{\beta-1}}
{(c-(a+t)^2)^{2\beta}}\\\times
   {}_2F_1\!\left[\begin{matrix}\alpha-\beta-1,-\beta\\
\alpha\end{matrix};\frac{(a+t)t}{c-a(a+t)}\right]\,
t^{\alpha-1}\,(1-t)^{\beta-1}\,dt.
\end{multline}
\end{Theorem}

We consider a few important special cases.
First, it is clear that Theorem~\ref{betat} reduces to the
classical beta integral evaluation \eqref{betaf} when either $c\to\infty$
or $a\to\infty$.
Other cases of interest concern the limits $c\to 0$ and $a\to 0$.
If $c\to 0$ \eqref{betat} reduces to
\begin{equation}\label{spec1}
\frac{\Gamma(\alpha)\Gamma(\beta)}{\Gamma(\alpha+\beta)}=
a^\beta(a+1)^{\beta+1}\int_0^1
\frac{t^{\alpha-1}(1-t)^{\beta-1}}{(a+t)^{2\beta+1}}\,
{}_2F_1\!\left[\begin{matrix}\alpha-\beta-1,-\beta\\
\alpha\end{matrix};\frac{-t}a\right]dt,
\end{equation}
which can easily
be recovered as a special case of Erd\'elyi's~\cite{Er} fractional
integral formula
(see also \cite[p.~112, Th.~2.9.1]{AAR})
\begin{multline}\label{erdelyi}
   {}_2F_1\!\left[\begin{matrix}a,b\\
c\end{matrix};x\right]=\frac{\Gamma(c)}{\Gamma(\mu)\Gamma(c-\mu)}\,
\int_0^1 t^{\mu-1}\,(1-t)^{c-\mu-1}\,(1-xt)^{\lambda-a-b}\\\times
{}_2F_1\!\left[\begin{matrix}\lambda-a,\lambda-b\\
\mu\end{matrix};xt\right]
{}_2F_1\!\left[\begin{matrix}a+b-\lambda,\lambda-\mu\\
c-\mu\end{matrix};\frac{(1-t)x}{1-xt}\right]\,dt.
\end{multline}
Indeed, \eqref{erdelyi} reduces to \eqref{spec1} when one
does the replacements $\lambda\mapsto\alpha$, $\mu\mapsto\alpha$,
$a\mapsto\beta+1$, $b\mapsto\alpha+\beta$, $c\mapsto\alpha+\beta$,
and $x\mapsto-1/a$, and then simplifies.
Some $q$-extensions of \eqref{erdelyi} and two other
fractional integral representations for hypergeometric functions
in \cite{Er} are derived in \cite{GG}.

For $a\to 0$, \eqref{betat} reduces to
\begin{multline}\label{spec5}
\frac{\Gamma(\alpha)\Gamma(\beta)}{\Gamma(\alpha+\beta)}=
(c-1)c^\beta\int_0^1
\frac{(c-t)^{\beta-1}}{(c-t^2)^{2\beta}}\,t^{\alpha-1}\,
(1-t)^{\beta-1}\\\times
{}_2F_1\!\left[\begin{matrix}\alpha-\beta-1,-\beta\\
\alpha\end{matrix};\frac{t^2}c\right]dt,
\end{multline}
which we were  unable to find in the literature.

Some special cases of \eqref{betat} where the ${}_2F_1$ in the
integrand can be simplified are $\alpha=\beta+1$, which is
\eqref{spec2}, and $\alpha=\beta$, which is \eqref{spec3}.

Next, we consider the beta-type integral evaluation
obtained from letting $q\to 1$ in \eqref{qbetat2}.
\begin{Theorem}\label{thbetat2}Let $\Re(\beta)>\max(0,m-1)$. Then
\begin{multline}\label{betat2}
\frac{\Gamma(\beta)\Gamma(\beta)}{2\,\Gamma(2\beta)}=
(c-(a+1)^2)
\int_0^1 \frac{(c-a(a+t))^\beta\,(c-(a+1)(a+t))^{\beta-1}}
{(c-(a+t)^2)^{2\beta}}\\\times
{}_2F_1\!\left[\begin{matrix}-\beta,-m\\
-2\beta\end{matrix};\frac{c-(a+t)^2}{(c-a(a+t))(1-et)}\right]
\left(\frac{1-et}{1-e}\right)^m
t^{\beta-m}\,(1-t)^{\beta-1}\,dt.
\end{multline}
\end{Theorem}

Note that \eqref{betat2} can be further rewritten
using Legendre's duplication formula
\begin{equation*}
\Gamma(2\beta)=\frac 1{\sqrt{\pi}}\,2^{2\beta-1}
\Gamma(\beta)\Gamma(\beta+\tfrac 12),
\end{equation*}
after which the left hand side becomes
\begin{equation*}
\frac{\sqrt{\pi}}{4^\beta}\,
\frac{\Gamma(\beta)}{\Gamma(\beta+\frac 12)}.
\end{equation*}

If $e=0$, \eqref{betat2} is equivalent to the $\alpha-\beta-1=-m$
case of \eqref{betat}, due to Pfaff's transformation
(cf.\ \cite[p.~79, Eq.~(2.3.14)]{AAR})
\begin{equation*}
{}_2F_1\!\left[\begin{matrix}-m,b\\
c\end{matrix};x\right]=\frac{(c-b)_m}{(c)_m}\,
{}_2F_1\!\left[\begin{matrix}-m,b\\
b+1-m-c\end{matrix};1-x\right].
\end{equation*}

If $c\to\infty$ or  $a\to\infty$, then \eqref{betat2} reduces
to the following beta-type integral evaluation:
\begin{equation}\label{spec4}
\frac{\Gamma(\beta)\Gamma(\beta)}{2\,\Gamma(2\beta)}=
\int_0^1 {}_2F_1\!\left[\begin{matrix}-\beta,-m\\
-2\beta\end{matrix};\frac 1{1-et}\right]
\left(\frac{1-et}{1-e}\right)^m
t^{\beta-m}\,(1-t)^{\beta-1}\,dt.
\end{equation}
For a conventional proof of \eqref{spec4},
expand the ${}_2F_1$ in powers of $1/(1-et)$, interchange
the order of summation and integration, and evaluate the
integrals using the $\lambda=\mu=a$ special case of \eqref{erdelyi},
interchange summations again, simplify by first using the
Gau{\ss} summation (cf.\ \cite[p.~66, Th.~2.2.2]{AAR}) and then
by using the binomial theorem. The details are left to the reader.

Other interesting cases of \eqref{betat2} are the specializations
$c=0$ or $a=0$.

Finally, observe that by performing various substitutions we may change
the form and path of integration of
the above integrals. In particular, using $t\mapsto
s/(s+1)$ these integrals then run over the half line $s\in[0,\infty)$.

\end{document}